\documentclass[a4paper,11pt,reqno]{article}

\usepackage[utf8]{inputenc}
\usepackage[T1]{fontenc}
\usepackage{lmodern}
\usepackage[english]{babel}
\usepackage{microtype}

\usepackage{amsmath,amssymb,amsfonts,amsthm}
\usepackage{mathtools,accents}
\usepackage{mathrsfs}
\usepackage{aliascnt}
\usepackage{braket}
\usepackage{bm}

\usepackage[a4paper,margin=3cm]{geometry}
\usepackage[citecolor=blue,colorlinks]{hyperref}

\usepackage{enumerate}
\usepackage{xcolor}

\makeatletter
\g@addto@macro\@floatboxreset\centering
\makeatother

\makeatletter
\def\newaliasedtheorem#1[#2]#3{
  \newaliascnt{#1@alt}{#2}
  \newtheorem{#1}[#1@alt]{#3}
  \expandafter\newcommand\csname #1@altname\endcsname{#3}
}
\makeatother

\numberwithin{equation}{section}

\newtheoremstyle{slanted}{\topsep}{\topsep}{\slshape}{}{\bfseries}{.}{.5em}{}

\theoremstyle{plain}
\newtheorem{theorem}{Theorem}[section]
\newaliasedtheorem{proposition}[theorem]{Proposition}
\newaliasedtheorem{lemma}[theorem]{Lemma}
\newaliasedtheorem{corollary}[theorem]{Corollary}
\newaliasedtheorem{counterexample}[theorem]{Counterexample}

\theoremstyle{definition}
\newaliasedtheorem{definition}[theorem]{Definition}
\newaliasedtheorem{question}[theorem]{Question}
\newaliasedtheorem{example}[theorem]{Example}
\newaliasedtheorem{conjecture}[theorem]{Conjecture}

\theoremstyle{remark}
\newaliasedtheorem{remark}[theorem]{Remark}

\let\altphi\phi
\let\phi\varphi
\let\varphi\altphi
\let\altphi\undefined

\newcommand{\loc}{{\rm loc}}

\newcommand{\Loc}{{\rm Loc}}
\newcommand{\spt}{{\rm spt}\,}

\DeclareMathOperator{\Lip}{Lip}
\DeclareMathOperator{\Lipb}{Lip_b}

\newfont{\tmpf}{cmsy10 scaled 2500}

\def\qedhere{\hfill $\Box$}

\begin{document}
\title{The pointed intrinsic flat distance \\ between locally integral current spaces}
\author{
Shu Takeuchi
\thanks{Tohoku University, \url{shu.takeuchi.q8@dc.tohoku.ac.jp}}} 
\maketitle

\begin{abstract}
In this note we define a distance 
between two pointed locally integral current spaces.
We prove that a sequence of pointed locally integral current spaces converges 
with respect to this distance
if and only if it converges in the sense of Lang-Wenger.
This enables us to state the compactness theorem by Lang-Wenger
for pointed locally integral current spaces
in terms of a distance function.
\end{abstract}

\tableofcontents

\section{Introduction}

Ambrosio and Kirchheim defined currents in metric spaces in \cite{AmbrosioKirchheim}.
We begin with a rough review of currents in metric spaces,
see the next section
for details.
Let $X$ be a complete separable metric space,
$k$ be a nonnegative integer and
${\cal D}_{\rm AK}^k(X):=\Lipb(X)\times(\Lip(X))^k$.
We call a multilinear function $T : {\cal D}_{\rm AK}^k(X)\to\mathbb{R}$ a $k$-dimensional current in $X$
if $T$ possesses
continuity, locality (see Definition \ref{metfunc}) and finiteness of the mass.
``Finiteness of the mass'' means that 
there exists a finite Borel measure $\mu$ on $X$ such that 
\begin{align}\label{mass-est}
|T(f , \pi_1 , \dots , \pi_k)|\leq\prod^k_{i=1}\Lip(\pi_i)\int_X|f|\,d\mu
\end{align}
holds for all $(f , \pi_1 , \dots , \pi_k)\in{\cal D}_{\rm AK}^k(X)$.
$\|T\|$ stands for the minimal measure $\mu$ satisfying (\ref{mass-est}),
and we define the mass of $T$ by
${\bf M}(T):=\|T\|(X)$. 
Let $T_1$ and $T_2$ be $k$-dimensional integral currents in $X$.
We define the flat distance between $T_1$ and $T_2$ in $X$ by
\begin{align*}
{\cal F}^X(T_1 , T_2):=\inf_{U , V}({\bf M}(U)+{\bf M}(V)),
\end{align*}
where the infimum is taken over all
$U\in{\bf I}_k(X)$ and $V\in{\bf I}_{k+1}(X)$ with $T_1-T_2=U+\partial V$
and ${\bf I}_k(X)$ denotes the set of all $k$-dimensional integral currents in $X$.


In order to discuss integral currents in different metric spaces,
Sormani and Wenger defined the intrinsic flat distance in \cite{SormaniWenger}.  
Roughly speaking,
the intrinsic flat distance between two $k$-dimensional integral current spaces $(X_1 , d_1 , T_1)$ and $(X_2 , d_2 , T_2)$
is defined by
\begin{align*}
d_{{\cal F}}((X_1 , d_1 , T_1) , (X_2 , d_2 , T_2)):=\inf{\cal F}^Z(\phi_{1\#}T_1 , \phi_{2\#}T_2)
\end{align*}
where the infimum is taken over all complete metric spaces $(Z , d)$
and isometric embeddings $\phi_1 : (X_1 , d_1)\hookrightarrow (Z , d)$ and $\phi_2 : (X_2 , d_2)\hookrightarrow (Z , d)$.
Furthermore, they proved a compactness theorem with respect to this distance.
The statement of the theorem is as follows: 
fix a sequence of complete separable metric spaces ${\{(X_n , d_n)\}}_n$ 
which are
uniformly bounded and uniformly totally bounded.
Let
$k\geq1$ and $T_n\in{\bf I}_k(X_n)$.
If we assume $\sup_n({\bf M}(T_n)+{\bf M}(\partial T_n))<\infty$,
then there exist an integral current space $(X , d , T)$ and a subsequence $n(j)$ such that
$(X_{n(j)} , d_{n(j)} , T_{n(j)})$ converges to $(X , d , T)$ with respect to the intrinsic flat distance.
Note that $(X , d)$ is not necessarily isometric to the Gromov-Hausdorff limit space of $(X_{n(j)} , d_{n(j)})$
(see Figure 2 in \cite{SormaniWenger}).
Also note that even if $(X_n , d_n , T_n)$ converges to $(X , d , T)$ with respect to $d_{{\cal F}}$,
the condition $\sup_n({\bf M}(T_n)+{\bf M}(\partial T_n))<\infty$ does not always hold.
It is not difficult to construct
a sequence ${\{T_n\}}_n\subset{\bf I}_k(\mathbb{R}^{k+1})$
which converges with respect to flat distance
and satisfies $\sup_n{\bf M}(T_n)=\infty$.

Now, 
let us discuss the case 
when the current spaces are
pointed locally integral current spaces,
which may have infinite mass. 
Let $X$ be a complete separable metric space, $x\in X$, 
$k\geq1$ and $T\in{\bf I}_{\Loc , k}(X)$,
where ${\bf I}_{\Loc , k}(X)$ denotes the set of all $k$-dimensional locally integral currents in $X$.
We call $(X , x , T)$ a $k$-dimensional pointed locally integral current space,
and ${\cal M}^k_*$ denotes the space of $k$-dimensional pointed locally integral current spaces.
In \cite{LangWenger}, 
Lang and Wenger studied the following convergence in ${\cal M}^k_*$:
for a sequence 
${\{(X_n , x_n , T_n)\}}_n\subset{\cal M}^k_*$,
we say that
$(X_n , x_n , T_n)$ converges to $(Z , z , T)\in{\cal M}^k_*$ 
if there is an isometric embedding $\phi_n : X_n\hookrightarrow Z$
such that $d_Z(\phi_n(x_n) , z)\to0 \ (n\to\infty)$
and $\phi_{n\#}T_n$ converges to $T$ in the local flat topology (Definition \ref{local-flat-top}).
Then they
proved
a compactness theorem for pointed locally integral current spaces 
with respect to this convergence.
It is a natural question that 
this convergence can be written in terms of a distance function
as in \cite{SormaniWenger},
that is, this convergence is metrizable or not.
By using the same idea as the intrinsic flat distance in \cite{SormaniWenger} 
and the pointed Gromov-Hausdorff distance in \cite{Gromovgroup}, 
we define the pointed intrinsic flat distance $d_{{\cal F}*}$ in ${\cal M}^k_*$,
whose convergence is compatible with the above convergence
(Definition \ref{pifd-def}, Proposition \ref{locflat--pifd} and Proposition \ref{pifd--locflat}):
consequently,  by the compactness theorem 
(Theorem 1.1 in \cite{LangWenger}: see \cite{LangWenger}, \cite{Wengercptness} , \cite{Wengerisoperimetric} 
and \cite{Wengerflat} for the proof), 
we obtain the following theorem:

\begin{theorem}\label{cptness}
Let $k\geq1$.
Assume that a sequence of $k$-dimensional pointed locally integral current spaces ${\{(X_n , x_n , T_n)\}}_n\subset{\cal M}^k_*$ satisfies
\begin{align*}
\sup_{n}(\|T_n\|+\|\partial T_n\|)(\bar{B}_r(x_n))<\infty
\end{align*}
for all $r>0$, where $\bar{B}_r(x_n)$ is the closed ball of radius $r$ centered at $x_n$.
Then there exist a subsequence ${\{(X_{n(j)} , x_{n(j)} , T_{n(j)})\}}_j$ and $(Z , z , T)\in{\cal M}^k_*$ such that
\begin{align*}
d_{{\cal F}*}((X_{n(j)} , x_{n(j)} , T_{n(j)}) , (Z , z , T))\to0
\end{align*}
as $j\to\infty$.
\qedhere\end{theorem}
In particular, 
we see that the convergence 
defined in \cite{LangWenger} for pointed locally integral current spaces
is ``intrinsic.''
Finally, we mention the $L^2_{\rm loc}$-convergence of orientations of Riemannian manifolds.
Let $(X , x , {\cal H}^k)$ be a Ricci limit space 
and ${\{(X^k_n , x_n , {\cal H}^k)\}}_n$ be a  sequence of $k$-dimensional Riemannian manifolds with
${\rm Ric}_{X_n}\geq-(k-1)$ and ${\cal H}^k(B_1(x_n))\geq v>0$.
Let $\omega_n$ (resp. $\omega$) be an orientation of $X_n$ (resp. $X$).
Then we see that
$(X_n , x_n , {\cal H}^k)$ mGH-converges to $(X , x , {\cal H}^k)$
and $\omega_n$ $L^2_{\rm loc}$-converges to $\omega$ simultaneously
if and only if $X_n$ converges to $X$ as pointed locally integral current spaces with respect to $d_{{\cal F}*}$ (see \cite{Honda} for details).

\section{Locally integral currents in metric spaces}
In this section, we recall the definition of $k$-dimensional locally integral currents and related notions.
See \cite{LangWenger} for details. 

Throughout this section, let $X$ be a complete separable metric space.
We define classes of Lipschitz functions as follows:
\begin{align*}
\Lip_{\Loc}(X)&:=\{f : X\to\mathbb{R} \ ; \ f|_B {\rm \ is \ Lipschitz \ for \ any \ bounded \ set \ }B\subset X\}, \\
\Lip(X)&:=\{f : X\to\mathbb{R} \ ; \ f {\rm \ is \ Lipschitz}\}, \\
\Lip_1(X)&:=\{f\in\Lip(X) \ ; \ |f(x)-f(y)|\leq d(x , y) {\rm \ for \ all \ } x , y\in X \}, \\
\Lipb(X)&:=\{f\in\Lip(X) \ ; \ f {\rm \ is \ bounded}\}, \\
\Lip_{\rm B}(X)&:=\{f\in\Lipb(X) \ ; \ \spt f{\rm \ is \ bounded}\},
\end{align*}
where $\spt f:=\overline{\{x\in X \ ; \ f(x)\neq0\}}$.
More generally, for a metric space $X'$,  we define
\begin{align*}
\Lip_{\Loc}(X , X')&:=\{f : X\to X' \ ; \ f|_B {\rm \ is \ Lipschitz \ for \ any \ bounded \ set \ }B\subset X\}, \\
\Lip(X , X')&:=\{f : X\to X' \ ; \ f {\rm \ is \ Lipschitz}\}.
\end{align*}
$\Lip(f)$ denotes the Lipschitz constant of $f\in\Lip(X)$, 
that is, $\Lip(f):=\sup\{|f(x)-f(y)|/d(x , y) \ ; \ x , y\in X , x\neq y\}$. 
For $r>0$, $B_r(x)$
(resp. $\bar{B}_r(x)$)
denotes the open
(resp. closed) 
ball centered at $x\in X$ with radius $r$.
Similarly, for a subset $A\subset X$, $B_r(A)$ (resp. $\bar{B}_r(A)$) denotes the 
open (resp. closed) $r$-neighborhood of $A$.

For $k\geq0$,
let ${\cal D}^k(X):=\Lip_{\rm B}(X)\times(\Lip_{\Loc}(X))^k$.
If $X^n$ is an $n$-dimensional Riemannian manifold, by Rademacher's theorem, 
any element $(f , \pi_1 , \dots , \pi_k)\in{\cal D}^k(X)$ determines a $k$-dimensional 
differential form $f\,d\pi_1\wedge\dots\wedge d\pi_k$ ${\cal H}^n$-a.e. on $X$.
Thus
we write $(f , \pi_1 , \dots , \pi_k)$
as $f\,d\pi_1\wedge\dots\wedge d\pi_k$ or $f\,d\pi$ for short.
\begin{definition}\label{metfunc}
A function $T : {\cal D}^k(X)\to\mathbb{R}$ is called a {\bf $k$-dimensional metric functional on $X$} 
if the following properties hold: \\
(i)(multilinearity) $T$ is multilinear. \\
(ii)(continuity) 
$\lim_{j\to\infty}T(f\,d\pi^j)=T(f\,d\pi)$ holds whenever
$\pi_i^j$ pointwisely converges to $\pi_i$ for any $i=1 , \dots , k$ with 
$\sup_j\Lip(\pi^j_i|_B)<\infty$ for any bounded set $B\subset X$. \\
(iii)(locality) If $\pi_i$ is constant on $\bar{B}_{\delta}(\spt f)$ for some $\delta>0$, then $T(f\,d\pi)=0$.
\qedhere\end{definition}
A typical example is as follows:
\begin{definition}\label{standardex}
Let $\theta\in L^1_\loc(\mathbb{R}^k)$, 
then 
a function $[\theta] : {\cal D}^k(\mathbb{R}^k)\to\mathbb{R}$
defined by
\begin{align*}
[\theta](f\,d\pi_1\wedge\dots\wedge d\pi_{k}):=\int_{\mathbb{R}^k}f\theta\det(\nabla\pi)\,d{\cal L}^k
\end{align*}
is a $k$-dimensional metric functional on $\mathbb{R}^k$,
where ${\cal L}^k$ denotes the Lebesgue measure on ${\mathbb{R}^k}$.
\qedhere\end{definition}
Now we define the pushforward, the restriction and the boundary of a metric functional.
\begin{definition}\label{operation}
Let $T$ be a $k$-dimensional metric functional on $X$. \\
(i) Let $\phi\in\Lip_{\Loc}(X , X')$,
and assume that for any bounded set $A\subset X'$ $\phi^{-1}(A)$ is also bounded.
Then we define a $k$-dimensional metric functional $\phi_{\#}T$ on $X'$ by
\begin{align*}
(\phi_{\#}T)(f , \pi_1 ,\dots , \pi_k):=T(f\circ\phi , \pi_1\circ\phi , \dots , \pi_k\circ\phi).
\end{align*}
$\phi_{\#}T$ is called the {\bf pushforward of $T$ with respect to $\phi$}. \\
(ii)
Let $l\in[0 , k]$ be an integer and 
$g\,d\tau
\in(\Lip_{\Loc}(X))^{l+1}$.
We define a $(k-l)$-dimensional metric functional $T\llcorner(g\,d\tau)$ on $X$ by
\begin{align*}
(T\llcorner(g\,d\tau))(f\,d\pi):=T(fg , \tau_1 , \dots , \tau_l , \pi_1 , \dots , \pi_{k-l}).
\end{align*}
$T\llcorner(g\,d\tau)$ is called the {\bf restriction of $T$ to $g\,d\tau$}. \\
(iii)
Let $k\geq1$.
We define a $(k-1)$-dimentional metric functional $\partial T$ by
\begin{align*}
(\partial T)(f\,d\pi):=T(1\,df\wedge d\pi_1\wedge\dots\wedge d\pi_{k-1}).
\end{align*}
$\partial T$ is called the {\bf boundary of $T$}.
\qedhere\end{definition}

Now we define the mass of a $k$-dimensional metric functional $T$.
\begin{definition}\label{metfunc-mass}
Let $T$ be a $k$-dimensional metric functional on $X$.
For any open set $O\subset X$, define
\begin{align*}
\|T\|(O):=\sup\sum^N_{j=1}T(f^j\,d\pi_1^j\wedge\dots\wedge d\pi_k^j)
\end{align*}
where  the supremum is taken over all $N\in\mathbb{N}$ and 
$f^j\,d\pi_1^j\wedge\dots\wedge d\pi_k^j\in\Lip_{\rm B}(X)\times(\Lip_1(X))^k$
with
$\spt f^j\subset O$ and $\sum^N_{j=1}|f^j|\leq1$. 
For all $A\subset X$, put
\begin{align*}
\|T\|(A):=\inf\{\|T\|(O) \ ; \ O\supset A {\rm \ is \ open}\}.
\end{align*} 
\qedhere\end{definition}
It is natural to ask whether $\|T\| : 2^X\to\mathbb{R}_{\geq0}\cup\{\infty\}$
is an outer measure or not. 
The following lemma gives an answer to this question.
See Proposition 2.2 in \cite{LangWenger} for the proof.
\begin{lemma}\label{outermeas}
Let $T$ be a $k$-dimensional metric functional on $X$.
Assume that for any bounded open set $O\subset X$ and $\epsilon>0$,
there exists a compact set $K\subset O$ such that $\|T\|(O\backslash K)<\epsilon$. 
Then $\|T\|$ is an outer measure satisfying that \\
{\rm (i)} any Borel set is $\|T\|$-measurable, \\
{\rm (ii)} for any $A\in 2^X$, there exists a Borel set $B$ such that $B\supset A$ and $\|T\|(A)=\|T\|(B)$, \\
{\rm (iii)} there exists a $\sigma$-compact set $\Sigma\subset X$
such that $\|T\|(\Sigma^c)=0$.
\qedhere\end{lemma}
For a metric functional $T$ satisfying the assumption of Lemma \ref{outermeas}, let
\begin{align*}
\spt T:=\{x\in X \ ; \ \|T\|(B_r(x))>0 {\rm \ for \ all \ } r>0\}.
\end{align*}
Using Definition \ref{metfunc-mass} and Lemma \ref{outermeas}, we define metric currents with locally finite mass.
Note that it may satisfy $\|T\|(X)=\infty$,  while \cite{SormaniWenger} deals with metric currents with $\|T\|(X)<\infty$.
\begin{definition}\label{def-current}
We say that a $k$-dimensional metric functional 
$T$ is a {\bf metric current with locally finite mass}
 if $\|T\|(O)<\infty$ holds for any bounded open set $O\subset X$
and the assumption of Lemma \ref{outermeas} holds.
${\bf M}_{\Loc , k}(X)$ denotes the vector space consisting of all $k$-dimensional metric currents with locally finite mass.
\qedhere\end{definition}
The next proposition will play a key role later,
see Proposition 2.3 in \cite{LangWenger} for the proof.
\begin{proposition}
Let $T\in{\bf M}_{\Loc , k}(X)$.
Then 
\begin{align}\label{currentestimate}
|T(f\,d\pi)|\leq\prod^k_{i=1}\Lip(\pi_i|_{\spt f})\int_X|f|\,d\|T\|
\end{align}
holds for all $f\,d\pi\in{\cal D}^k(X)$.
\qedhere\end{proposition}
In order to extend the domain of $T\in{\bf M}_{\Loc , k}(X)$,
let us use the following notation:
\begin{align*}
\mathscr{B}^\infty_{\Loc}(X)&:=\left\{f : X\to\mathbb{R} \ ; \ {\rm Borel \ measurable}, \sup_A|f|<\infty {\rm \ for \ any \ bounded \ set \ }
A\subset X\right\}, \\
\mathscr{B}^\infty(X)&:=\left\{f : X\to\mathbb{R} \ ; \ {\rm Borel \ measurable}, \sup_X|f|<\infty\right\}, \\
\mathscr{B}^\infty_{\rm B}(X)&:=\left\{f\in\mathscr{B}^\infty(X) \ ; \ \spt f {\rm \ is \ bounded}\right\}.
\end{align*}
In the following, we use the next lemma.
\begin{lemma}\label{Lip-app}
Let $T\in{\bf M}_{\Loc , k}(X)$, $f\in\mathscr{B}^\infty_{\rm B}(X)$ and $N$ be a bounded neighborhood of $\spt f$.
Then there exists ${\{f_n\}}_n\subset \Lip_{\rm B}(X)$ such that $\spt f_n\subset N$ 
holds for all $n$ and $f_n\to f$ in $L^1(\|T\|)$.
\qedhere\end{lemma}
\begin{proof}
Fix an arbitrary $n\in\mathbb{N}$.
Let $R>0$ be a positive number such that $B_{R}(\spt f)\subset N$. 
By Definition \ref{def-current},
there exists a compact set $K\subset B_{R}(\spt f)$ such that $\|T\|(B_{R}(\spt f)\backslash K)<(n(1+\sup_X|f|))^{-1}$.
Then we have
\begin{align*}
\int_X|f-f\chi_K|\,d\|T\|=\int_{B_R(\spt f)\backslash K}|f|\,d\|T\|<\frac{1}{n},
\end{align*}
where $\chi_K$ is the indicator function of $K$.
Moreover, by the definition of the integral, there exist finite Borel sets $B_1 , \dots , B_m$ and real numbers
$a_1 , \dots  , a_m$ such that  $B_i\subset K$ for all $i=1 , \dots , m$ and that
\begin{align*}
\int_X|f\chi_K-\sum^m_{i=1}a_i\chi_{B_i}|\,d\|T\|<\frac{1}{n}.
\end{align*}
Finally, for all $i$, there exists $f_{n , i}\in\Lip_{\rm B}(X)$ such that
\begin{align*}
\int_X|a_i\chi_{B_i}-f_{n , i}|\,d\|T\|<\frac{1}{nm}.
\end{align*}
Indeed, for any $\epsilon>0$, 
define $f_{n , i , \epsilon}\in\Lip_{\rm B}(X)$ by $f_{n , i , \epsilon}(x):=a_i\max\{0 , 1-\epsilon^{-1}d(x , B_i)\}$.
Then one can find such functions by letting $\epsilon\to0$.
Taking $f_n:=\sum^m_{i=1}f_{n , i}$ completes the proof.
\end{proof}
\begin{definition}\label{Borel-ext}
Let $T\in{\bf M}_{\Loc , k}(X)$. 
For $(f , \pi_1 , \dots , \pi_k)\in\mathscr{B}^\infty_{\rm B}(X)\times(\Lip_{\Loc}(X))^k$,
define
\begin{align}\label{borellimit}
T(f , \pi_1 , \dots , \pi_k):=\lim_{n\to\infty}T(f_n\,d\pi)
\end{align}
where $N$ is a bounded neighborhood of $\spt f$, and ${\{f_n\}}_n\subset\Lip_{\rm B}(X)$ satisfies $\spt f_n\subset N$ and $f_n\to f$
in $L^1(\|T\|)$, 
as in Lemma \ref{Lip-app}.
\qedhere\end{definition}
The limit in (\ref{borellimit}) exists since  ${\{T(f_n\,d\pi)\}}_n$ is a Cauchy sequence by (\ref{currentestimate}) and Lemma \ref{Lip-app}. 
If we take another bounded neighborhood $N'$ of $\spt f$ and another sequence ${\{g_n\}}_n\subset \Lip_{\rm B}(X)$
with $\spt g_n\subset N'$ and $g_n\to f$ in $L^1(\|T\|)$, then
\begin{align*}
|T(f_n\,d\pi)-T(g_n\,d\pi)|\leq\prod^k_{i=1}\Lip(\pi_i|_{N\cup N'})\int_X|f_n-g_n|\,d\|T\|\to0 \ (n\to\infty)
\end{align*}
holds.
Therefore
the limit in (\ref{borellimit})
does not depend on the choice of $N$ and ${\{f_n\}}_n$. 

Now we 
define the restriction to $g\,d\tau\in\mathscr{B}^\infty_{\Loc}(X)\times(\Lip_{\Loc}(X))^l$:
\begin{definition}
Let $T\in{\bf M}_{\Loc , k}(X)$, $l\in[0 , k]$ be an integer and 
$g\,d\tau\in\mathscr{B}^\infty_{\Loc}(X)\times(\Lip_{\Loc}(X))^l$.
We define $T\llcorner(g\,d\tau)\in{\bf M}_{\Loc , k-l}(X)$ by
\begin{align}\label{borelextension}
(T\llcorner(g\,d\tau))(f\,d\pi):=T(fg , \tau_1 , \dots , \tau_l , \pi_1 , \dots , \pi_{k-l})
\end{align}
for $f\,d\pi\in{\cal D}^{k-l}(X)$,
where the right-hand side of (\ref{borelextension}) is well-defined by Definition \ref{Borel-ext}.
$T\llcorner(g\,d\tau)$ is called the {\bf restriction of $T$ to $g\,d\tau$}.
\qedhere\end{definition}
\begin{definition}
Let $T\in{\bf M}_{\Loc , k}(X)$. For a Borel set $A\subset X$, 
we define the {\bf restriction of $T$ to $A$} by
\begin{align*}
T\llcorner\,A:=T\llcorner\,\chi_A
\end{align*}
where $\chi_A$ is the indicator function of $A$.
\qedhere\end{definition}
Using above notions, we define $k$-dimensional integral currents.
We say $S\subset X$ is a {\bf compact $k$-rectifiable set}
if there exist finite compact sets $K_1 , \dots , K_N$ in $\mathbb{R}^k$ and $\pi_i\in\Lip(K_i , X) \ (i=1 , \dots , N)$ such that
$S=\bigcup^N_{i=1}\pi_i(K_i)$.
\begin{definition}
A $k$-dimensional metric functional $T$ is said to be a {\bf $k$-dimensional locally integer rectifiable current}
if the following two conditions are satisfied: \\
(i) For any bounded open set $O\subset X$, 
we see that $\|T\|(O)<\infty$ and that for any $\epsilon>0$ there exists a compact $k$-rectifiable set $K$ such that 
$\|T\|(O\backslash K)<\epsilon$. \\
(ii) For any bounded Borel set $B\subset X$ and $\pi\in\Lip(X , \mathbb{R}^k)$, there exists $\theta\in L^1(\mathbb{R}^k , \mathbb{Z})$
such that $\pi_{\#}(T\llcorner\,B)=[\theta]$.

${\cal I}_{\Loc , k}(X)$ denotes the set of all $k$-dimensional locally integer rectifiable currents.
\qedhere\end{definition}
\begin{definition}
$T\in{\cal I}_{\Loc , k}(X)$ is said to be a {\bf $k$-dimensional locally integral current}
if $\partial T\in{\cal I}_{\Loc , k-1}(X)$ holds.
${\bf I}_{\Loc , k}(X)$ denotes the set of all $k$-dimensional locally integral currents.
\qedhere\end{definition}
Finally, we introduce a notion of convergence in ${\bf I}_{\Loc , k}(X)$.
\begin{definition}\label{local-flat-top}
Let $T_n , T\in{\bf I}_{\Loc , k}(X)$.
We say {\bf $T_n$ converges to $T$ in the local flat topology}
if for any bounded closed set $B\subset X$ there exists $U_n\in{\bf I}_{\Loc , k}(X)$ and 
$V_n\in{\bf I}_{\Loc , k+1}(X)$ such that $T_n-T=U_n+\partial V_n$ and $(\|U_n\|+\|V_n\|)(B)\to 0 \ (n\to\infty)$.
\qedhere\end{definition}

\section{The pointed intrinsic flat distance}

In this section,
we introduce the pointed intrinsic flat distance $d_{{\cal F}*}$(Definition \ref{pifd-def}).
It is a distance between pointed locally integral current spaces,
which may have infinite mass.
 Note that the intrinsic flat distance $d_{\cal F}$ in \cite{SormaniWenger} deals with integral current spaces,
which have finite mass.

Let $k\in\mathbb{Z}_{\geq1}$. 
For a complete separable metric space $X$, $x\in X$ and $T\in {\bf I}_{\Loc ,k}(X)$, we call a triplet $(X , x , T)$
a $k$-dimensional pointed locally integral current space.
The set of all $k$-dimensional pointed locally integral current spaces is denoted by ${\cal M}^k_*$.
\begin{definition}\label{pifd-def}
Let $(X_1 , x_1 , T_1) , (X_2 , x_2 , T_2)\in{\cal M}^k_*$.
We define the pointed intrinsic flat distance between $(X_1 , x_1 , T_1)$ and $(X_2 , x_2 , T_2)$ by
\begin{align*}
d_{{\cal F}*}((X_1 , x_1 , T_1) , (X_2 , x_2 , T_2)):=\min\left\{\widetilde{d_{{\cal F}*}}((X_1 , x_1 , T_1) , (X_2 , x_2 , T_2)) , \frac{1}{2}\right\}
\end{align*}
where $\widetilde{d_{{\cal F}*}}((X_1 , x_1 , T_1) , (X_2 , x_2 , T_2))$ is the infimum of $\epsilon>0$ satisfying following conditions:
there exist a complete metric space $Z$ and an isometric embedding $\phi_i : X_i\hookrightarrow Z \ (i=1 , 2)$ such that
\begin{itemize}
\item[(i)] $d_Z(\phi_1(x_1) , \phi_2(x_2))<\epsilon$,
\item[(ii)] for $i=1 , 2$, there exist $U_i\in{\bf I}_{\Loc , k}(Z)$ and $V_i\in{\bf I}_{\Loc , k+1}(Z)$ 
such that $\phi_{1\#}T_1-\phi_{2\#}T_2=U_i+\partial V_i$ and $(\|U_i\|+\|V_i\|)(\bar{B}_{1/\epsilon}(\phi_i(x_i)))<\epsilon$.
\end{itemize}
\qedhere\end{definition}
Let us check that $d_{{\cal F}*}$ is a pseudodistance on ${\cal M}^k_*$.
We recall the gluing of two metric spaces along same isometric images.
\begin{lemma}\label{two-glue}
Let $X , Z^1 , Z^2$ be metric spaces and $\phi^i : X\hookrightarrow Z^i \ (i=1 , 2)$ be an isometric embedding. 
Define $d : (Z^1\sqcup Z^2)\times(Z^1\sqcup Z^2)\to[0 , \infty)$ by
\begin{align*}
d(z , z'):=
\begin{cases}
d_{Z^1}(z , z') &(z , z'\in Z^1), \\
d_{Z^2}(z , z') &(z , z'\in Z^2), \\
\inf_{x\in X}(d_{Z^1}(z , \phi^1(x))+d_{Z^2}(\phi^2(x) , z')) &(z\in Z^1 , z'\in Z^2), \\
\inf_{x\in X}(d_{Z^2}(z , \phi^2(x))+d_{Z^1}(\phi^1(x) , z')) &(z\in Z^2 , z'\in Z^1),
\end{cases}
\end{align*}
then $d$ is a pseudodistance on $Z^1\sqcup Z^2$.
Moreover, let $(Z^1\sqcup Z^2)/d$ be the quotient metric space,
that is, let $(Z^1\sqcup Z^2)/d$ be the quotient space with respect to the equivalence relation defined by
\begin{align*}
z\sim z'\Leftrightarrow d(z , z')=0 \ (z , z'\in Z^1\sqcup Z^2).
\end{align*}
Then the canonical inclusion 
$\iota^i : Z^i\hookrightarrow (Z^1\sqcup Z^2)/d \ (i=1 , 2)$ is an isometric embedding.
In the following, $Z^1\sqcup_XZ^2$ denotes $(Z^1\sqcup Z^2)/d$.
\qedhere\end{lemma}
\begin{proof}
It is enough to check that $d$ satisfies the triangle inequality;
\begin{align}\label{two-glue-met-tri}
d(z , z'')\leq d(z , z')+d(z' , z'')  \ {\rm for \ any} \ z , z' , z''\in Z^1\sqcup Z^2.
\end{align}
If $z\in Z^1, z' , z''\in Z^2$, (\ref{two-glue-met-tri}) follows from that for all $x\in X$ 
\begin{align*}
d(z , z'')
\leq d_{Z^1}(z , \phi^1(x))+d_{Z^2}(\phi^2(x) , z'') 
\leq d_{Z^1}(z , \phi^1(x))+d_{Z^2}(\phi^2(x) , z')+d_{Z^2}(z' , z'')
\end{align*}
holds.
If $z , z''\in Z^1 , z'\in Z^2$, (\ref{two-glue-met-tri}) follows from that for all $x , x'\in X$
\begin{align*}
d(z , z'')
&\leq d_{Z^1}(z , \phi^1(x))+d_{Z^1}(\phi^1(x) , \phi^1(x'))+d_{Z^1}(\phi^1(x') , z'') \\
&=d_{Z^1}(z , \phi^1(x))+d_{Z^2}(\phi^2(x) , \phi^2(x'))+d_{Z^1}(\phi^1(x') , z'') \\
&=d_{Z^1}(z , \phi^1(x))+d_{Z^2}(\phi^2(x) , z')+d_{Z^2}(z' , \phi^2(x'))+d_{Z^1}(\phi^1(x') , z'')
\end{align*}
holds. 
Similarly, we can prove (\ref{two-glue-met-tri}) in the remaining cases.
\end{proof}
\begin{proposition}
$d_{{\cal F}*}$ is a pseudodistance on ${\cal M}^k_*$.
\qedhere\end{proposition}
\begin{proof}
We use a simplified notation $d_{{\cal F}*}(X_i , X_j)$ (resp. $\widetilde{d_{{\cal F}*}}(X_i , X_j)$) 
instead of $d_{{\cal F}*}((X_i , \\ x_i , T_i) , (X_j , x_j , T_j))$ (resp. $\widetilde{d_{{\cal F}*}}((X_i , x_i , T_i) , (X_j , x_j , T_j))$).
It is enough to check that $d_{{\cal F}*}$ satisfies the triangle inequality.
For $(X_i , x_i , T_i)\in{\cal M}^k_* \ (i=1 , 2 , 3)$, we have to show that
\begin{align*}
d_{{\cal F}*}(X_1 , X_3)\leq d_{{\cal F}*}(X_1 , X_2)+d_{{\cal F}*}(X_2 , X_3)
\end{align*}
holds.
Without loss of generality, we can assume that $\widetilde{d_{{\cal F}*}}(X_1 , X_2)<1/2$
and $\widetilde{d_{{\cal F}*}}(X_2 , X_3) \\ <1/2$.
For any sufficiently small $\delta>0$ with $\widetilde{d_{{\cal F}*}}(X_1 , X_2)+\delta<1/2$
and $\widetilde{d_{{\cal F}*}}(X_2 , X_3)+\delta<1/2$,
there exist $\epsilon_1<\widetilde{d_{{\cal F}*}}(X_1 , X_2)+\delta$ and
$\epsilon_2<\widetilde{d_{{\cal F}*}}(X_2 , X_3)+\delta$ such that the following holds: \\
(I) there exist a complete metric space $Z^1$ and an isometric embedding $\phi^1_i : X_i\hookrightarrow Z^1 \ (i=1 , 2)$ such that
\begin{itemize}
\item[(I-i)] $d_{Z^1}(\phi^1_1(x_1) , \phi^1_2(x_2))<\epsilon_1$,
\item[(I-ii)] for $i=1 , 2$, there exist $U^1_i\in{\bf I}_{\Loc , k}(Z^1)$ and $V^1_i\in{\bf I}_{\Loc , k+1}(Z^1)$ 
such that $\phi^1_{1\#}T_1-\phi^1_{2\#}T_2=U^1_i+\partial V^1_i$ 
and $(\|U^1_i\|+\|V^1_i\|)(\bar{B}_{1/\epsilon_1}(\phi^1_i(x_i)))<\epsilon_1$.
\end{itemize} 
(I\hspace{-.1em}I) there exist a complete metric space $Z^2$ and 
an isometric embedding $\phi^2_i : X_i\hookrightarrow Z^2 \ (i=2 , 3)$ such that
\begin{itemize}
\item[(I\hspace{-.1em}I-i)] $d_{Z^2}(\phi^2_2(x_2) , \phi^2_3(x_3))<\epsilon_2$,
\item[(I\hspace{-.1em}I-ii)] for $i=2 , 3$, there exist $U^2_i\in{\bf I}_{\Loc , k}(Z^2)$ and $V^2_i\in{\bf I}_{\Loc , k+1}(Z^2)$ 
such that $\phi^2_{2\#}T_2-\phi^2_{3\#}T_3=U^2_i+\partial V^2_i$ 
and $(\|U^2_i\|+\|V^2_i\|)(\bar{B}_{1/\epsilon_2}(\phi^2_i(x_i)))<\epsilon_2$.
\end{itemize}
Now, we apply Lemma \ref{two-glue} to $\phi^1_2 : X_2\hookrightarrow Z^1$ and $\phi^2_2 : X_2\hookrightarrow Z^2$.
Let $Z:=Z^1\sqcup_{X_2}Z^2$ and $\iota^i : Z^i\hookrightarrow Z \ (i=1 , 2)$ be the canonical inclusion, then we define 
an isometric embedding $\phi_i : X_i\hookrightarrow Z \ (i=1 , 2 , 3)$ by
\begin{align*}
\phi_1&:=\iota^1\circ\phi^1_1, \\
\phi_2&:=\iota^1\circ\phi^1_2(=\iota^2\circ\phi^2_2), \\
\phi_3&:=\iota^2\circ\phi^2_3.
\end{align*}
Then, 
\begin{align*}
d_Z(\phi_1(x_1) , \phi_3(x_3))
&\leq d_Z(\phi_1(x_1) , \phi_2(x_2))+d_Z(\phi_2(x_2) , \phi_3(x_3)) \\
&=d_{Z^1}(\phi^1_1(x_1) , \phi^1_2(x_2))+d_{Z^2}(\phi^2_2(x_2) , \phi^2_3(x_3)) \\
&<\epsilon_1+\epsilon_2
\end{align*}
and for $i=1 , 2$,
\begin{align*}
\phi_{1\#}T_1-\phi_{3\#}T_3
&=(\phi_{1\#}T_1-\phi_{2\#}T_2)+(\phi_{2\#}T_2-\phi_{3\#}T_3) \\
&=\iota^1_{\#}(\phi^1_{1\#}T_1-\phi^1_{2\#}T_2)+\iota^2_{\#}(\phi^2_{2\#}T_2-\phi^2_{3\#}T_3) \\
&=\iota^1_{\#}(U^1_i+\partial V^1_i)+\iota^2_{\#}(U^2_{i+1}+\partial V^2_{i+1}) \\
&=(\iota^1_{\#}U^1_i+\iota^2_{\#}U^2_{i+1})+\partial(\iota^1_{\#}V^1_i+\iota^2_{\#}V^2_{i+1})
\end{align*}
holds.
Moreover, we have
\begin{align*}
&(\|\iota^1_{\#}U^1_1+\iota^2_{\#}U^2_2\|+\|\iota^1_{\#}V^1_1+\iota^2_{\#}V^2_2\|)(\bar{B}_{1/(\epsilon_1+\epsilon_2)}(\phi_1(x_1))) \\
&\leq(\|\iota^1_{\#}U^1_1\|+\|\iota^1_{\#}V^1_1\|)(\bar{B}_{1/(\epsilon_1+\epsilon_2)}(\phi_1(x_1)))
+(\|\iota^2_{\#}U^2_2\|+\|\iota^2_{\#}V^2_2\|)(\bar{B}_{1/(\epsilon_1+\epsilon_2)}(\phi_1(x_1))) \\
&\leq(\|\iota^1_{\#}U^1_1\|+\|\iota^1_{\#}V^1_1\|)(\bar{B}_{1/\epsilon_1}(\phi_1(x_1)))
+(\|\iota^2_{\#}U^2_2\|+\|\iota^2_{\#}V^2_2\|)(\bar{B}_{1/\epsilon_2}(\phi_2(x_2))) \\
&\leq(\|U^1_1\|+\|V^1_1\|)(\bar{B}_{1/\epsilon_1}(\phi^1_1(x_1)))
+(\|U^2_2\|+\|V^2_2\|)(\bar{B}_{1/\epsilon_2}(\phi^2_2(x_2))) \\
&<\epsilon_1+\epsilon_2,
\end{align*}
where we used the fact that $\bar{B}_{1/(\epsilon_1+\epsilon_2)}(\phi_1(x_1))\subset\bar{B}_{1/\epsilon_2}(\phi_2(x_2))$,
which follows from inequalities; $d_Z(\phi_1(x_1) , \phi_2(x_2))<\epsilon_1$, $\epsilon_1<1/2$ and $\epsilon_2<1/2$.
Similarly,  we have
\begin{align*}
(\|\iota^1_{\#}U^1_2+\iota^2_{\#}U^2_3\|+\|\iota^1_{\#}V^1_2+\iota^2_{\#}V^2_3\|)(\bar{B}_{1/(\epsilon_1+\epsilon_2)}(\phi_3(x_3)))
<\epsilon_1+\epsilon_2.
\end{align*}

From the above, we have
\begin{align*}
d_{{\cal F}*}(X_1 , X_3)\leq\epsilon_1+\epsilon_2<d_{{\cal F}*}(X_1 , X_2)+d_{{\cal F}*}(X_2 , X_3)+2\delta.
\end{align*}
Letting $\delta\to0$ completes the proof.
\end{proof}
Let us check that 
the convergence in Theorem 1.1 in \cite{LangWenger}
coincides with 
that by the pointed intrinsic flat distance.
First, the following proposition shows that the convergence in \cite{LangWenger} implies that of $d_{{\cal F}*}$.
\begin{proposition}\label{locflat--pifd}
Let ${\{(X_n , x_n , T_n)\}}_n\subset{\cal M}^k_*$ be a sequence of $k$-dimensional pointed locally integral current spaces.
Assume that there exist $(Z , z , T)\in{\cal M}^k_*$ and an isometric embedding $\phi_n : X_n\hookrightarrow Z$ such that
$\phi_n(x_n)\to z \ (n\to\infty)$ and $\phi_{n\#}T_n$ converges to $T$ in the local flat topology.
Then $d_{{\cal F}*}((X_n , x_n , T_n) , (Z , z , T))\to0 \ (n\to\infty)$.
\qedhere\end{proposition}
\begin{proof}
Fix an arbitrary $\epsilon\in(0 , 1/2)$. Then for sufficiently large $n$, $d_{Z}(\phi_n(x_n) , z)<\epsilon$ holds, 
and there exist $U_n\in{\bf I}_{\Loc , k}(Z)$ and $V_n\in{\bf I}_{\Loc , k+1}(Z)$ such that 
$\phi_{n\#}T_n-T=U_n+\partial V_n$ and
$(\|U_n\|+\|V_n\|)(\bar{B}_{1+1/\epsilon}(z))<\epsilon$ (recall Definition \ref{local-flat-top}).
Since $\bar{B}_{1/\epsilon}(\phi_n(x_n))\subset \bar{B}_{1+1/\epsilon}(z)$, the conclusion follows.
\end{proof}
Next, let us check that the converse of Proposition \ref{locflat--pifd} is also true.
We use the following lemma in order to prove that (Proposition \ref{pifd--locflat}).
\begin{lemma}\label{sigma-glue-met}
Let $Z , Z^i \ (i=1 , 2 , \dots)$ be complete metric spaces and $\phi^i : Z\hookrightarrow Z^i \ (i=1 , 2 , \dots)$ be an isometric embedding.
We define $d : (\bigsqcup_{i=1}^\infty Z^i)\times(\bigsqcup_{i=1}^\infty Z^i)\to[0 , \infty)$ by
\begin{align*}
d(z , z'):=
\begin{cases}
d_{Z^i}(z , z') &(i=j), \\
\inf_{\bar{z}\in Z}(d_{Z^i}(z , \phi^i(\bar{z}))+d_{Z^j}(\phi^j(\bar{z}) , z')) &(i\neq j),
\end{cases}
\end{align*}
where $z\in Z^i , z'\in Z^j$.
Then $d$ is a pseudodistance on $\bigsqcup_{i=1}^\infty Z^i$.
Moreover, $(\bigsqcup_{i=1}^\infty Z^i)/d$ is a complete metric space 
and the canonical inclusion $\iota^i : Z^i\hookrightarrow (\bigsqcup_{i=1}^\infty Z^i)/d$
is an isometric embedding. 
In the following, 
we will use a simplified notation $(\bigsqcup_{i=1}^\infty Z^i)/Z$ 
instead of $(\bigsqcup_{i=1}^\infty Z^i)/d$.
\qedhere\end{lemma}
\begin{proof}
Let us check that $d$ satisfies the triangle inequality. For $z\in Z^i , z'\in Z^j$ and $z''\in Z^k$, 
we have to show
\begin{align}\label{sigma-glue-met-tri}
d(z , z'')\leq d(z , z')+d(z' , z'').
\end{align}
If $i , j$ and $k$ are different each other, then (\ref{sigma-glue-met-tri}) holds because we have
\begin{align*}
d(z , z'')
&\leq d_{Z^i}(z , \phi^i(\bar{z}))+d_{Z^k}(\phi^k(\bar{z}) , z'') \\
&\leq d_{Z^i}(z , \phi^i(\bar{z}))+d_{Z^k}(\phi^k(\bar{z}) , \phi^k(\tilde{z}))+d_{Z^k}(\phi^k(\tilde{z}) , z'') \\
&=d_{Z^i}(z , \phi^i(\bar{z}))+d_{Z^j}(\phi^j(\bar{z}) , \phi^j(\tilde{z}))+d_{Z^k}(\phi^k(\tilde{z}) , z'') \\
&\leq d_{Z^i}(z , \phi^i(\bar{z}))+d_{Z^j}(\phi^j(\bar{z}) , z')+d_{Z^j}(z' , \phi^j(\tilde{z}))+d_{Z^k}(\phi^k(\tilde{z}) , z'')
\end{align*}
for any $\bar{z} , \tilde{z}\in Z$.
Otherwise, 
(\ref{sigma-glue-met-tri}) holds by Lemma \ref{two-glue}.

Next, we show that $(\bigsqcup_{i=1}^\infty Z^i)/d$  is a complete metric space.
Let ${\{z_n\}}_n\subset(\bigsqcup_{i=1}^\infty Z^i)/d$ be a Cauchy sequence. 
If $\#\{z_n ; z_n\in Z^i\}=\infty$ for some $i$,
we can choose a subsequence ${\{z_{n(j)}\}}_j$ such that $z_{n(j)}\in Z^i$ holds for all $j$.
Since ${\{z_{n(j)}\}}_j$ converges to some $z\in Z^i$, ${\{z_n\}}_n$ is a convergent sequence in $(\bigsqcup_{i=1}^\infty Z^i)/d$.
If $\#\{z_n ; z_n\in Z^i\}<\infty$ for any $i$, 
after taking a subsequence, 
there exists a
strictly increasing sequence ${\{i(n)\}}_n$ such that $z_n\in Z^{i(n)}$.
Since ${\{z_n\}}_n$ is a Cauchy sequence, we have
\begin{align*}
d_{Z^{i(n)}}(z_n , \phi^{i(n)}(Z))\to0 \ (n\to\infty).
\end{align*}
Let $\bar{z}_n\in Z$ with $d_{Z^{i(n)}}(z_n , \phi^{i(n)}(\bar{z}_n))<d_{Z^{i(n)}}(z_n , \phi^{i(n)}(Z))+1/2^n$, then
\begin{align*}
d_Z(\bar{z}_n , \bar{z}_m)
&=d(\phi^{i(n)}(\bar{z}_n) , \phi^{i(m)}(\bar{z}_m)) \\
&\leq d(\phi^{i(n)}(\bar{z}_n) , z_n)+d(z_n , z_m)+d(z_m , \phi^{i(m)}(\bar{z}_m)) \\
&=d_{Z^{i(n)}}(\phi^{i(n)}(\bar{z}_n) , z_n)+d(z_n , z_m)+d_{Z^{i(m)}}(z_m , \phi^{i(m)}(\bar{z}_m)) \\
&\to0 \ (n , m\to\infty).
\end{align*}
Hence ${\{\bar{z}_n\}}_n\subset Z$ is a Cauchy sequence and converges to some $\bar{z}\in Z$. 
Here,
\begin{align*}
d(\phi^{i(n)}(\bar{z}) , z_n)
&=d_{Z^{i(n)}}(\phi^{i(n)}(\bar{z}) , z_n) \\
&\leq d_{Z^{i(n)}}(\phi^{i(n)}(\bar{z}) , \phi^{i(n)}(\bar{z}_n))+d_{Z^{i(n)}}(\phi^{i(n)}(\bar{z}_n) , z_n) \\
&=d_Z(\bar{z} , \bar{z}_n)+d_{Z^{i(n)}}(\phi^{i(n)}(\bar{z}_n) , z_n) \\
&\to0 \ (n\to\infty)
\end{align*}
holds.
Since $d(\phi^{i(n)}(\bar{z}) , \phi^{i(m)}(\bar{z}))=0$ for any $n , m\in\mathbb{N}$, the conclusion follows.
\end{proof}
Now we show the converse of Proposition \ref{locflat--pifd}.
\begin{proposition}\label{pifd--locflat}
Assume that $(X_n , x_n , T_n) , (Z , z , T)\in{\cal M}^k_*$ satisfy $d_{{\cal F}*}((X_n , x_n , T_n) , (Z , \\ z , T))\to0 \ (n\to\infty)$.
Then there exist $(Z' , z' , T')\in{\cal M}^k_*$ and an isometric embedding $\phi_n : X_n\hookrightarrow Z'$ such that
$\phi_n(x_n)\to z' \ (n\to\infty)$ and that $\phi_{n\#}T_n$ converges to $T'$ in the local flat topology.
Moreover, we can take $(Z' , z' , T')$ in such a way that
there exists an isometric embedding $\phi : Z\hookrightarrow Z'$
satisfying $\phi(z)=z'$ and $\phi_{\#}T=T'$.
\qedhere\end{proposition}
\begin{proof}
Let $\epsilon_n:=\widetilde{d_{{\cal F}*}}((X_n , x_n , T_n) , (Z , z , T))$. By Definition \ref{pifd-def},
there exists $\delta_n\in[\epsilon_n , \epsilon_n+1/2^n)$ 
(if $\epsilon_n=0$, $\delta_n\in(0 , 1/2^n)$) 
satisfying following conditions: \\
There exist a complete metric space $Z^n$ and  isometric embeddings $\phi^n : X_n\hookrightarrow Z^n$,
$\psi^n : Z\hookrightarrow Z^n$ such that 
\begin{itemize}
\item[(i)] $d_{Z^n}(\phi^n(x_n) , \psi^n(z))<\delta_n$,
\item[(ii)] there exist 
$U_n , \tilde{U}_n\in{\bf I}_{\Loc , k}(Z^n)$, $V_n , \tilde{V}_n\in{\bf I}_{\Loc , k+1}(Z^n)$ such that 
$\phi^n_{\#}T_n-\psi^n_{\#}T=U_n+\partial V_n=\tilde{U}_n+\partial \tilde{V}_n$,
that $(\|U_n\|+\|V_n\|)(\bar{B}_{1/\delta_n}(\phi^n(x_n)))<\delta_n$ and that
$(\|\tilde{U}_n\|+\|\tilde{V}_n\|)(\bar{B}_{1/\delta_n}(\psi^n(z)))<\delta_n$.
\end{itemize}

Now we apply Lemma \ref{sigma-glue-met} to $\psi^n : Z\hookrightarrow Z^n$.
Let $Z':=(\bigsqcup^\infty_{n=1}Z^n)/Z$ and $\varphi_n : Z^n\hookrightarrow Z' , \phi : Z\hookrightarrow Z'$
be the canonical isometric embeddings.
If we take $z':=\phi(z)$, $T':=\phi_{\#}T$ and $\phi_n:=\varphi_n\circ\phi^n$, the conclusion follows.
In fact, 
\begin{align*}
d_{Z'}(z' , \phi_n(x_n))=d_{Z^n}(\psi^n(z) , \phi^n(x_n))<\delta_n
\end{align*}
implies the convergence of the reference point. 
Moreover, since
\begin{align*}
\phi_{n\#}T_n-\phi_{\#}T=\varphi_{n\#}(\phi^n_{\#}T_n-\psi^n_{\#}T)=\varphi_{n\#}(\tilde{U}_n+\partial\tilde{V}_n)
=(\varphi_{n\#}\tilde{U}_n)+\partial(\varphi_{n\#}\tilde{V}_n),
\end{align*}
$\phi_{n\#}T_n$ converges to $T'$ in the local flat topology because
\begin{align*}
(\|\varphi_{n\#}\tilde{U}_n\|+\|\varphi_{n\#}\tilde{V}_n\|)(B)
&\leq(\|\varphi_{n\#}\tilde{U}_n\|+\|\varphi_{n\#}\tilde{V}_n\|)(\bar{B}_{1/\delta_n}(z')) \\
&\leq(\|\tilde{U}_n\|+\|\tilde{V}_n\|)(\bar{B}_{1/\delta_n}(\psi_n(z))) \\
&<\delta_n
\end{align*}
where $B\subset Z'$ is an arbitrary bounded closed set and $n$ is a sufficiently large number satisfying 
$B\subset \bar{B}_{1/\delta_n}(z')$.
\end{proof}
In the end of our note, 
we discuss the case of $d_{{\cal F}*}=0$.
\begin{proposition}
If $d_{{\cal F}*}((X , x , T) , (X' , x' , T'))=0$,
then there exists an isometry $\psi : \{x\}\cup\spt T\to\{x'\}\cup\spt T'$ such that $\psi(x)=x'$ and $\psi_{\#}T=T'$.
\qedhere\end{proposition}
\begin{proof}
Since $d_{{\cal F}*}((X , x , T) , (X' , x' , T'))=0$, there exists a
sequence of $\delta_n>0$ with $\lim_{n\to\infty}\delta_n=0$ and 
following conditions:
there exist a complete metric space $Z_n$ and isometric embeddings 
$\phi_n : X\hookrightarrow Z_n , \phi'_n : X'\hookrightarrow Z_n$ such that
\begin{itemize}
\item[(i)] $d_{Z_n}(\phi_n(x) , \phi'_n(x'))<\delta_n$,
\item[(ii)] there exist $U_n , \tilde{U}_n\in{\bf I}_{\Loc , k}(Z_n)$ and $V_n , \tilde{V}_n\in{\bf I}_{\Loc , k+1}(Z_n)$ 
such that $\phi_{n\#}T-\phi'_{n\#}T'=U_n+\partial V_n=\tilde{U}_n+\partial\tilde{V}_n$ , 
$(\|U_n\|+\|V_n\|)(\bar{B}_{1/\delta_n}(\phi_n(x)))<\delta_n$
and $(\|\tilde{U}_n\|+\|\tilde{V}_n\|)(\bar{B}_{1/\delta_n} \\ (\phi'_n(x')))<\delta_n$.
\end{itemize}

Now we apply Lemma \ref{sigma-glue-met} to $\phi'_n : X'\hookrightarrow Z_n$ and let $Z:=(\bigsqcup^\infty_{n=1}Z_n)/X'$.
Then let $\iota_n : Z_n\hookrightarrow Z$ and $\psi' : X'\hookrightarrow Z$ be canonical isometric embeddings, 
and one can check $\psi_n(x)\to\psi'(x') \ (n\to\infty)$ and $\psi_{n\#}T$ converges to $\psi'_{\#}T'$ in the local flat topology
where $\psi_n:=\iota_n\circ\phi_n$.
Now the conclusion follows from Proposition 1.1 in \cite{LangWenger}.
\end{proof}

In particular, $d_{{\cal F}*}$ is a distance function on the quotient space ${\cal M}^k_*/\sim$,
where $\sim$ is an equivalence relation on ${\cal M}^k_*$ defined as follows:
we say $(X , x , T)\sim(X' , x' , T')$ if and only if
 there exists an isometry $\psi : \{x\}\cup\spt T\to\{x'\}\cup\spt T'$ such that $\psi(x)=x'$ and $\psi_{\#}T=T'$.

\end{document}